\documentclass[twoside,leqno,fleqn]{article}
\usepackage{amsmath}
\usepackage{amssymb}

\usepackage{calligra}
\usepackage[T1]{fontenc}

\usepackage[shortlabels]{enumitem}
\usepackage[onehalfspacing]{setspace}

\newcommand{\AUTOR}{C. Ga\ss ner}
\newcommand{\TITEL}{AC and the Independence of WO in Second-Order Henkin Logic, Part I}
\def\S{\Sigma}
\def\I{{\cal I}}
\def\H{\mathfrak{H}}
\def\G{\mathfrak{G}}
\def\F{\mbox{\scriptsize\calligra T\!\!\!\!\!\hspace{0.005cm}F\!\!\!\!\!T\!\!\!\!\!\hspace{0.02cm}F}\,\,\,}

\newcommand{\bbbn}{\mathbb{N}}

\newcommand{\mbm}[1]{\mbox{\boldmath{$#1$}}}
\newcommand{\mbmss}[1]{\mbox{\scriptsize\boldmath{$#1$}}}
\newcommand{\mbmty}[1]{\mbox{\tiny\boldmath{$#1$}}}

\newtheorem{satz}{Satz}[section] 
\newtheorem{lemma}[satz]{Lemma}
\newtheorem{proposition}[satz]{Proposition}
\newtheorem{corollary}[satz]{Corollary}
\newtheorem{theorem}[satz]{Theorem}

\newtheorem{defi}[satz]{Definition}

\newtheorem{agree}[satz]{Agreement}

\begin{document}

\newcounter{li}

\thispagestyle{empty}
\begin{center} {\Large\bf AC and the Independence of WO \vspace{0.4cm}\\ in Second-Order Henkin Logic, Part I}\vspace{0.6cm}\\{\bf Christine Ga\ss ner}\footnote{I thank Michael Rathjen for the discussion on various questions related to the paper of Siskind, Mancosu, and Shapiro (2020) in Leeds. Moreover, I thank Arno Pauly and Vasco Brattka for helpful technical advice on the scope of pre-publications of proofs and the organizers of the Colloquium Logicum  in Konstanz in 2022.} {\vspace{0.2cm}\\University of Greifswald, Germany, 2024\\ gassnerc@uni-greifswald.de}\\\end{center} 

\begin{abstract} This article is concerned with the Axiom of Choice (AC) and the well-ordering theorem (WO) in second-order predicate logic with Henkin interpretation (HPL).  We consider a principle of choice introduced by Wilhelm Ackermann (1935) and discussed also by David Hilbert and Ackermann (1938), by  Günter Asser (1981), and by Benjamin Siskind, Paolo Mancosu, and Stewart  Shapiro (2020).  
The discussion is restricted to so-called Henkin-Asser structures of second order. The language used is a many-sorted first-order language with identity.  In particular, we give some of the technical details for a proof of the independence of WO from the so-called Ackermann axioms in HPL presented at the Colloquium Logicum in 2022.
\end{abstract}

\section{Introduction}\label{SectionEinf}
In \cite{1c}, G\"unter Asser gave a detailed introduction to a {\sf Henkin interpretation} of second-order predicate logic (denoted by HPL) and to a class of so-called {\sf Henkin structures} that are first-order characterizable. In this framework, individuals and predicates can form a domain of discourse, a so-called universe of a many-sorted first-order structure. The individuals are of sort 0 and the sort of a predicate is its arity. The language used is a many-sorted first-order language with identity and with variables for individuals and variables for predicates. Such a language makes it possible, for example, to define principles of choice. Let us remark that the discussed model-theoretic semantics is a specific form of Henkin interpretations which were detailed described by Leon Henkin (1950), Stewart Shapiro (1991), Jouko Väänänen (2019), and others. The possible interpretations of a formula can lead to different theories. Our language of HPL is a language with an identity symbol. In other cases, new formulas with the identity symbol are introduced as abbreviations for formulas without an identity symbol. It is therefore possible that the identity relations on predicates are not considered as components of the models. In Section 7 of \cite{1c}, Asser discussed the possibility of constructing Henkin structures using permutations.

Motivated by \cite{1c}, the dissertation \cite{Gass84} continues this approach and investigates the relative strength of second-order formulas, whose first-order counterparts --- such as the well-ordering theorem and the trichotomy law --- are known as equivalent to the Axiom of Choice (denoted by AC)  in ZF, and other modifications of AC in {\rm HPL} (for an overview see \cite{Gass94}).
On the one hand, \cite{1c} and \cite{Gass84} show that certain proofs of known relationships and many independence proofs can be transferred from
ZFA to {\rm HPL}. 
On the other hand, other proofs using various techniques to form higher-rank sets from lower-rank sets cannot be transferred from ZFA to HPL.
Here, we recall the necessary basic notions and describe the Fraenkel-Mostowski-Specker-Asser method for constructing Henkin structures of second order only shortly. Then, we deal with questions about certain formulations of the principle of choice that go back to Ackermann and Hilbert \cite{{HiAc38}} and Asser \cite{1c} and that are also discussed by Siskind, Mancosu, and Shapiro in \cite{Siskind} and their strength in {\rm HPL} and we also discuss some proofs. 

\section {Second-Order Predicate Logic}\label{Abschn_Abbr} 

We recall some important notions and statements. Our  language and our structures are of signature $\mbm{\sigma}^{(2)}$.

 \fbox{\parbox{11.3cm}{
\sf \small
 $ \mbm{\sigma}^{(2)}=(\bbbn,\emptyset,\emptyset,\{\mbm{\sf r}_n\mid n\geq 1\})$  with 
\begin{itemize}
\item sorts $0,1,2,\ldots$
\item a set of sorts: $ \bbbn=\{0,1,2,\ldots\}$
\item {\em sorts} $ \mbm{\sf r}_n=(n,0,0,\ldots,0)\in \bbbn^{n+1}$ ($n\geq 1$)
\end{itemize}
}}
\subsection{The language}
For the basic notions and prerequisites see \cite{1c,Gass84, Gass94}. Here, we give only a short overview. The used symbols and variables are the following.

\fbox{\parbox{11.3cm}{
\sf \small
\begin{itemize}
\item{\em a symbol for logic with identity}: $= $
\item {\em individual variables} of sort $0$, for elements of sort $0$: $ x_0,x_1,\ldots$,
\item {\em predicate variables} of sort $n\geq 1$,  for $n$-ary predicates: $ A_0^n, A_1^n,\ldots$ 
\item an $(n+1)$-ary {\em relation symbol} $ \mbm{\sf R}_n$ of sort $\mbm{\sf r}_n$ ($n\geq 1$)
\end{itemize}
}}

\vspace{0.3cm}
Let ${\cal L}^{(2)} $ be the set of second-order PL-II formulas defined as follows.

\fbox{\parbox{11.3cm}{
 \sf \small

\begin{list}{(\arabic{li})}
{\usecounter{li} \labelwidth0.5cm \leftmargin1.3cm \itemsep2pt plus1pt
\topsep1pt plus1ptminus1pt
\labelsep0.4cm\parsep0.5pt plus0.1pt minus0.1pt \itemindent-0.2cm}
\item $\mbm{\sf R}_n(A_i^n,x_{i_1},\ldots, x_{i_n}) $ and $x_i=x_j$ and $A_i^n=A_j^n$ are formulas. 

\item For all formulas $H$, $H_1$, and $H_2$ 

\qquad $\neg H$, $ (H_1 \land H_2)$, $ (H_1 \lor H_2)$, $ (H_1 \to H_2)$, and $(H_1\leftrightarrow H_2)$

are formulas.

\item For each formula $H(x_i)$, in which $x_i$ occurs only free, 

\qquad $ \forall {x_i}H(x_i)$ and $ \exists {x_i}H(x_i)$

are formulas.

\item For each formula $H(A^n_i)$, in which $A^n_i$ occurs only free, 

\qquad $ \forall {A^n_i}H(A^n_i)$ and $ \exists {A^n_i}H(A^n_i)$ 

are formulas.

\item $H$ is a formula if and only if it is it because of (1), \ldots, (4).
\end{list}

\vspace{0.1cm}
Later, we use $A_i^nx_{i_1}\cdots x_{i_n} $ instead of $\mbm{\sf R}_n(A_i^n,x_{i_1},\ldots, x_{i_n}) $.

\hfill ($n\geq 1$ and $i, j,i_1,\ldots, i_n\geq 0$)
}}

\vspace{0.3cm}
 Each atomic formula $\mbm{\sf R}_n(A_{i}^{n},x_{i_1}, \ldots, x_{i_n})$ in ${\cal L}^{(2)}$ can only refer to a uniquely determined relation and a predicate variable  $A_{i}^{n}$ of fixed arity $n \geq 1$. For simplicity, we often write $A_{i}^{n}x_{i_1}\cdots x_{i_n}$ instead of $\mbm{\sf R}_n(A_{i}^{n},x_{i_1}, \ldots, x_{i_n})$ for any $n\geq 1$ and all $i,i_1,\ldots,i_n\geq 0$. 
Moreover, we will also use other variables such as $x, x_i, x_i^{(j)},\ldots$, $y,y_i, y_i^{(j)},\ldots $ as individual variables (where $i,j, \ldots$ are placeholders for integers $\geq 0$), the (meta)\-variables $A,B,C,D, R,S,T$ as placeholders for predicate variables, the metavariables $H,H_1,\ldots$, $F,F_1,\ldots$, and $G,G_1,\ldots$ as placeholders for formulas and subformulas, notations such as $lo(T,A)$, $wo(T,A)$, and $choice_*^{n,m}(H)$ as abbreviations of certain formulas, and so on. 
In the following, $\mbm{x},\mbm{x}_i, \ldots$ in general are $n$-tuples given by $\mbm{x}=(x_{1}, \ldots,x_n)$, $\mbm{x}_i=(x_{1}^{(i)}, \ldots,x_{n}^{(i)})$, and so on. In formulas, we use $\mbm{x}$ as abbreviation for the list $x_1,\ldots,x_n$ and for the string $ x_1\cdots x_n$. This means that the string $A_i^n\mbm{x}$ stands for $A_i^n x_1 \cdots x_n $ and thus for $\mbm{\sf R}_n(A_{i}^{n},x_{1}, \ldots, x_{n})$ and $H(x_1,\ldots,x_n)$ can be abbreviated by $H(\mbm{x})$. Note that we write also $H(x_1, \ldots, x_n) $ or $H(\mbm{x}) $ instead of $H$ if $H$ is a second-order formula in which the variables $x_1, \ldots, x_n$ may occur and each occurrence of each of the variables $x_1, \ldots, x_n$ is free, and so on. Let ${\cal L}^{(2)}_{\mbmss{x},A_{i_1}^{n_1},\ldots,A_{i_s}^{n_s}}$ be the set of all formulas $H(\mbm{x},A_{i_1}^{n_1},\ldots,A_{i_s}^{n_s})$ in ${\cal L}^{(2)}$ in which $x_{1},\ldots,x_{n}$ and $A_{i_1}^{n_1},\ldots,A_{i_s}^{n_s}$ may occur only free. The sets ${\cal L}^{(2)}_{x_{1},\ldots,x_n}$, ${\cal L}^{(2)}_{A_{i_1}^{n_1},\ldots,A_{i_s}^{n_s}}$, and the like are defined by analogy. $\mbm{y},\mbm{y}_j, \ldots$ are $m$-tuples and $\mbm{y}$ can be used as an abbreviation for $y_1,\ldots,y_m$ and for $ y_1\cdots y_m$. $A_i^{n+m}\mbm{x}\mbm{y} $ stands for $A_i^{n+m} x_1\cdots x_ny_1\cdots y_m $ and $H(A_i^n,\mbm{x},\mbm{y})$ stands for $H(A_i^n,x_1,\ldots,x_n, y_1,\ldots,y_m)$, and so on. Moreover, let $( \mbm{x}\,.\,\mbm{y})$ be the $(n+m)$-tuple $(x_1,\ldots,x_n, y_1,\ldots,y_m)$. $\forall \mbm{x}$ means $\forall x_1\cdots \forall x_n$ and $\exists \mbm{x}$ means $\exists x_1\cdots \exists x_n$. Moreover, $\mbm{x}_1=\mbm{x}_2$ stands for $x_{1}^{(1)}=x_{1}^{(2)}\land \cdots \land x_{n}^{(1)}=x_{n}^{(2)}$. The abbreviation $\exists !! \mbm{x} H(\mbm{x})$ stands for \[\mbox{$\exists \mbm{x} H (\mbm{x}) \,\,\land\,\, \forall \mbm{x}_1\forall \mbm{x}_2 (H(\mbm{x}_1)\,\,\land\,\, H(\mbm{x}_2)\to \mbm{x}_1=\mbm{x}_2) $}.\] Note that other authors also use $\exists ^1$, $\exists ^{=1}$ or $\exists !$ instead of $\exists !!$. For a finite set $M$, let $\bigvee_{i\in M}H_i$ be the formula $H_{m_1}\lor \cdots\lor H_{m_s}$ and $\bigwedge_{i\in M}H_i$ be the formula $ H_{m_1}\land \cdots\land H_{m_s}$ if $M=\{m_1,\ldots, m_s\}$. Moreover, let $H\land \bigwedge_{i\in \emptyset}H_i=_{\rm df} H$, let $H\land \bigwedge_{i\in \{l,\ldots,k\}}H_i=_{\rm df} H$ if $k<l$, and the like.

 \subsection{Predicate structures}\label{DefPredicateStructures}

Each predicate structure considered here will be a many-sorted first-order structure of signature $\mbm{\sigma}^{(2)}$ that can be completely determined by its universe $ J_0\cup J_1\cup J_2\cup\cdots $ where $J_{0}$ is a non-empty domain of individuals and each non-empty domain $J_n$ ($n\geq 1$) contains $n$-ary predicates. For determining the relationships between the individuals and the predicates, we will define $(n+1)$-ary relations $\mbm{r}^{J_0}_{J_n}$ of sort $\mbm{\sf r}_n$ for all $n\geq 1$ in the same way. Each $\mbm{r}^{J_0}_{J_n}$ is a subset of $J_n\times J_0\times\cdots \times J_0$ and, for a fixed structure, we also write $\mbm{r}_n$ instead of $\mbm{r}^{J_0}_{J_n}$.
Our second-order predicate structures  $(\bigcup_{n\geq 0}J_n; \emptyset;\emptyset;(\mbm{r}_n)_{n\geq 1})$ can be described as follows. This means that we have the following.

\fbox{\parbox{11.3cm}{
\sf \small

\begin{itemize}
\item many-sorted first-order structures $\S$: $(\bigcup_{n\geq 0}J_n; \emptyset;\emptyset;(\mbm{r}_n)_{n\geq 1})$ 
\item signature: $\mbm{\sigma}^{(2)}$
\item elements $\xi_0,\xi_1,\ldots$ of sort $0$ in $J_0$
\item $n$-ary predicates $\alpha_0,\alpha_1,\ldots$ in $J_n$
\item relation $ \mbm{r}_n$ or  $ \mbm{r}_{J_n}^{J_0}$ (for interpreting the symbol $\mbm{\sf R}_n$)
\[(\alpha,\xi_1,\ldots,\xi_n)\in \mbm{r}_n \mbox{ if and only if } \alpha(\xi_1,\ldots,\xi_n)=true \mbox{ for } \alpha \in J_n\]
\[(\alpha,\xi_1,\ldots,\xi_n)\in \mbm{r}_n \mbox{ if and only if } (\xi_1,\ldots,\xi_n)\in \widetilde \alpha \hspace*{0.6cm}\mbox{ for } \alpha \in J_n\]
\end{itemize}
}}

\vspace{0.3cm}

For  any non-empty set $I$ and $n\geq 1$,  let $\mbm{\xi}$ and $\mbm{\xi}_i$ be abbreviations for the $n$-tuples $(\xi_1,\ldots,\xi_n)\in I^n$ and $(\xi_{i,1},\ldots,\xi_{i,n})\in I^n$, respectively, and the lists $\xi_1,\ldots,\xi_n$ and $\xi_{i,1},\ldots,\xi_{i,n}$, respectively, of individuals in $ I$. Let, for any fixed $m\geq 1$, $\mbm{\eta}$ be an abbreviation for an $m$-tuples $(\eta_1,\ldots,\eta_m)\in I^m$ of individuals $\eta_1,\ldots,\eta_m\in I$. In particular, we use the abbreviation $(\mbm{\xi}\,.\,\mbm{\eta})$ for the $(n+m)$-tuple $(\xi_1,\ldots,\xi_n,\eta_1,\ldots,\eta_m)\in I^{n+m}$ and so on.

\vspace{0.3cm}

Any $n$-ary function $\alpha: I^n\!\to \{true, false\}$ that assigns to each tuple $\mbm{\xi}\in I^n$ one of the values $true$ or $false$ is called an {\em$n$-ary predicate} ({\em on $I$}). 
Consequently, every predicate $\alpha: I^n\to \{true, false\}$ defines a relation $\widetilde \alpha=\{ \mbm{\xi}\in I^n\mid \alpha(\mbm{\xi})=true\}$. 
Let ${\rm pred}_n(I)$ be the domain of all $n$-ary predicates $\alpha: I^n\to \{true, false\}$ on $I$. Then, ${\rm pred}_n(I)$ is the set $\{\chi^I_{\mbmss{r}}\mid \mbm{r}\subseteq I^{n}\}$ with $\chi^I_{\mbmss{r}}(\mbm{\xi})= true$ if $\mbm{\xi}\in \mbm{r}$ and otherwise $\chi^I_{\mbmss{r}}(\mbm{\xi})= false$. Moreover, let ${\rm pred}(I)$ be the domain of all predicates on $I$. This means that ${\rm pred}(I)=\bigcup_{n\geq 1} {\rm pred}_n(I)$. Note that this union was denoted by $I^{\rm pred}$ in \cite{Gass94}.

Let $\mbm{r}^I_n$ be the set $\{(\alpha,\xi_1,\ldots,\xi_n) \mid \alpha \in {\rm pred}_n(I) \,\,\&\,\,\mbm{\xi}\in I^n\,\,\&\,\,\alpha(\mbm{\xi})=true \}$. Thus, $\mbm{r}^I_n $ is an $(n+1)$-ary relation of sort $\mbm{\sf r}_n$.
Let $J_{0}$ be an arbitrary non-empty set $I$ and, for any $n\geq 1$, let $J_n$ be a non-empty subset of ${\rm pred}_n(J_0)$. Let $\mbm{r}^{J_0}_{J_n}$ be defined by $\mbox{$\mbm{r}^{J_0}_{J_n}=\{(\alpha,\xi_1,\ldots,\xi_n)\in \mbm{r}^{J_0}_n \mid \alpha \in J_n\,\,\&\,\,\mbm{\xi}\in J_0^n \}$}$. Then, the structure $(\bigcup_{n\geq 0}J_n;\emptyset;\emptyset; (\mbm{r}^{J_0}_{J_n})_{ n\geq 1})$ is called a {\em predicate structure} ({\em over $J_0$}) and denoted by $(J_n)_{n\geq 0}$. 
Consequently, $\mbm{r}^{J_0}_{J_n}\subseteq \mbm{r}^I_n$ for all $n\geq 1$. The elements of $J_{0}$ are of {\em sort $0$}. The predicates in each {\em domain $J_n$} are of {\em sort $n$}. Then, $(\bigcup_{n\geq 0}J_n;\emptyset;\emptyset; (\mbm{r}^{J_0}_{J_n})_{ n\geq 1})$ is a structure of signature $\mbm{\sigma}^{(2)}$. It is a structure with the many-sorted set $\bigcup_{n\geq 0}J_n$ as {\em universe}, without constants, without functions, and with one $(n+1)$-ary relation $\mbm{r}^{J_0}_{J_n}$ of {\em sort $\mbm{\sf r}_n$} for any $n\geq 1$. 
 The complete predicate structure $(J_{n})_{n\geq 0}$ given by $J_0=I$ and $J_n= {\rm pred}_n(I)$ for all $n\geq 1$ is {\em the standard structure $\S_{\sf 1}^{I}$ over $I$}. Consequently,  for this many-sorted first-order structure, we have  $\S_{\sf 1}^{I}= (I\cup{\rm pred}(I);\emptyset;\emptyset; (\mbm{r}^I_n)_{ n\geq 1})$ and the structure $\widetilde{\S_{\sf 1}^{I}}$ given by $\widetilde{\S_{\sf 1}^{I}}=(I;\emptyset;\emptyset;(I^n) _{ n\geq 1})$ is the corresponding one-sorted first-order structure. Let ${\sf struc}^{({\rm m})}_{\rm pred}( I)$ be the class of all predicate structures $(J_n)_{n\geq 0}$ with $J_0=I$.

\subsection{Assignments and models}
Any predicate structure $\S$ allows us to evaluate all well-formed formulas of PL\,II by means of assignments $(f_n)_{n\geq 0}$ {\em in $\S$} (cf.\,\,\cite{1c}).

\begin{defi}[Assignments in ${\rm assgn}(\S)$] Let $\S$ be a structure $(J_n)_{n\geq 0}$ in ${\sf struc}_{\rm pred}^{({\rm m})}(J_0)$. Any {\em $\mbm{\sigma}^{(2)}$-assignment function} $f$ in $\S$ (shortly, {\em assignment} in $\S$) assigns to each variable $x_i$ an individual $\xi_i\in J_0$ and, for every $n\geq 1$, it assigns to each predicate variable $A^n_i$ an $n$-ary predicate $\alpha_i\in J_n$. Let ${\rm assgn}(\S)$ be the set of all $\mbm{\sigma}^{(2)}$-assignments in $\S$. 
\end{defi} 
$f$ determines a sequence $(f_n)_{n\geq 0}$ of functions $f_n$ with $f_0:\{x_i\mid i\geq 0\}\to J_0$ and $f_n:\{A^n_i\mid i\geq 0\}\to J_n$ ($n> 0$). Let $f_0( x_i)=f(x_i)$ for each variable $x_i$ ($i\geq 0$) and $f_n(A^n_i) =f(A^n_i)$ for each variable $A^n_i$ ($n\geq 1, i\geq 0$). Let $f$, $\bar f$, $\tilde f$, $f'$, and the like be names for assignments that determine sequences $( f_n)_{n\geq 0}$, $(\bar f_n)_{n\geq 0}$, $(\tilde f_n)_{n\geq 0}$, $(f'_n)_{n\geq 0}$ and so on, respectively. Let $f\big\langle{x\atop \xi} \big\rangle$ be an assignment like $f$, expect that it assigns to the individual variable $x$ the individual $\xi \in J_0$. This means that, for any assignments $f$ and $f'$ with $f'=f\big\langle{x\atop \xi} \big\rangle$, $f$ and $f'$ assign the same values to all other variables. Supposing $f'=f\big\langle{x\atop \xi} \big\rangle$, we have $f'(x)= \xi$, $f'(y)=f(y)$ if $x\not=y$, and $f'(A_i^n)=f(A_i^n)$ for all $n\geq 1$ and $i\geq 0$. For the tuple $\mbm{x}=(x_1,\ldots,x_n)$ and the tuple $\mbm{\xi}=(\xi_1,\ldots,\xi_n)$ in $J_0^n$, we use the notation $f\big\langle{\mbmss{x}\atop\mbmss{\xi}} \big\rangle$ for an assignment in order to express that this assignment assigns to $x_i$ in $ \mbm{x}$ the component $\xi_i$ in $ \mbm{\xi}$ ($i\leq n$). Let $f\big\langle{\mbmss{x}\atop\mbmss{\xi}} \big\rangle= f\big\langle{x_{1}\atop \xi_1} \big\rangle$ if $n=1$ and, for $n\geq 2$, let $ f\big\langle{x_{1}\atop \xi_1}{\cdots\atop \cdots}{x_{n}\atop \xi_n} \big\rangle=\big(f\big\langle{x_{1}\atop \xi_1}{\cdots\atop \cdots}{x_{n-1}\atop \xi_{n-1}} \big\rangle\big)\big\langle{x_{n}\atop \xi_n} \big\rangle$. Moreover, let $f\big\langle{\mbmss{x}\atop\mbmss{\xi}} \big\rangle= f\big\langle{x_{1}\atop \xi_1}{\cdots\atop \cdots}{x_{n}\atop \xi_n} \big\rangle$, $f\big\langle{A_i^n\atop \alpha}{\mbmss{x} \atop \mbmss{\xi}}\big\rangle=\big(f\big\langle {A_i^n\atop \alpha}\big\rangle\big)\big\langle{\mbmss{x}\atop\mbmss{\xi}}\big\rangle$, and so on. If $f'=f\big\langle{\mbmss{x}\atop \mbmss{\xi}} \big\rangle$ holds for some $f$, then we can also write $f'(\mbm{x})=\mbm{\xi}$ and the like. 

\begin{defi}[The interpretation $I_\S$] Let $\S$ be any predicate structure \linebreak $(J_n)_{n\geq 0}$. Now, let  $n\geq 1$. By $\S$,
the symbol $\mbm{\sf R}_n$ is interpreted as the relation $\mbm{r}^{J_0}_{J_n}\subseteq J_n\times J_0\times\cdots\times J_0$. To express this fact, we define $I_\S(\mbm{\sf R}_n)=\mbm{r}^{J_0}_{J_n}$ and write shortly $\S(\mbm{\sf R}_n)=\mbm{r}^{J_0}_{J_n}$.
\end{defi}

\begin{defi}[The valuation functions $\S_f$] Let $\S$ be a structure $(J_n)_{n\geq 0}$ in ${\sf struc}_{\rm pred}^{({\rm m})}(J_0)$ and $f$ be in ${\rm assgn}(\S)$. For any $H$ in ${\cal L}^{(2)}$, let $\S_f(H)$ --- denoted by $Wert_{\S}(H,f)$ in \cite{1c} and below also by $\S[H,f]$ --- be the {\em truth value of $H$ in $\S$ under $f$}. All truth values $\S_f(H)$ result from a recursive definition of the valuation function $\S_f:{\cal L}^{(2)}\to \{true, false\}$ assigning truth values to all {\rm PL-II} formulas $H$.
\end{defi}

\begin{defi}[The truth values of the atomic formulas in $\S$]For any predicate variables $A_i^n$ and $A_j^n$ of sort $n\geq 1 $ ($ i,j\geq 0$)  and any   individual variables $x_i, x_j, x_{i_1}, \ldots, x_{i_n}$  of sort 0 ($i,j,i_1,\ldots,i_n\geq 0$), let the truth values of the formulas $x_i=x_j$, $A_i^n =A_j^n$, and $\mbm{\sf R}_n(A_{i}^{n},x_{i_1}, \ldots, x_{i_n})$ of depth 0 in ${\cal L}^{(2)}$ be defined as follows. 
 \[\begin{array}{ll}
\vspace{0.2cm} \S_f(x_i=x_j)&=\,true \mbox{ iff } \,\, f(x_i)=f(x_j).\\
\vspace{0.2cm} \S_f(A_i^n=A_j^n)&=\,true \mbox{ iff } \,\,f(A_i^n)=f(A_j^n).\\ 
 \S_{f}(\mbm{\sf R}_n(A_{i}^{n},x_{i_1}, \ldots, x_{i_n}))\!\!\!&=\,true \mbox{ iff } \,\,(f_n(A^n_i ),f_0(x_{i_1}),\ldots,f_0(x_{i_n}))\in \S(\mbm{\sf R}_n).\end{array}\] 
\end{defi}
Note, that {\em iff} stands here for {\em if and only if}. Let $\in$ stand for ``{\em belongs to}''. $\in$ is only a relation in the metatheory, it does not belong to $\S$.

Because of $\S(\mbm{\sf R}_n)=\mbm{r}^{J_0}_{J_n}$, this means that we get $ \S_{f}(\mbm{\sf R}_n(A^n_i, x_{i_1}, \ldots, x_{i_n}))=f_n(A^n_i )(f_0(x_{i_1}),\ldots,f_0(x_{i_n}) )$ and, thus, \[\S_f(A^n_i x_{i_1} \cdots x_{i_n} ) =f_n(A^n_i )(f_0(x_{i_1}),\ldots,f_0(x_{i_n}) ).\] 
 $f_n(A^n_i )(f_0(x_{i_1}),\ldots,f_0(x_{i_n}) )$ stands for $(f_n(A^n_i ))(f_0(x_{i_1}),\ldots,f_0(x_{i_n}) )$ and consequently, for any $\mbm{\xi} \in J_0^n$ and any $\alpha\in J_n$, we get the equation \[\S_{f\langle{A_i^n\atop \alpha}{{x_{i_1}\atop \xi_1} \cdots {x_{i_n}\atop \xi_n}}\rangle}(A^n_i x_{i_1} \cdots x_{i_n}) =\alpha(\mbm{\xi}).\] 

\noindent Second, for all composite formulas $H$ in ${\cal L}^{(2)}$, the values $\S_f(H)$ can be defined --- by the principle of extensionality --- as is usual in second-order logic.

\begin{agree}\label{AgreeFreeVar} 
Let $x_1,\ldots,x_n$ and $A_{i_1}^{n_1}, \ldots,A_{i_r}^{n_r}$ be the only variables that may occur free in $H\in {\cal L}^{(2)}$. The relevant assignments can be specified directly and, thus, the symbol $f$ in $f\big\langle {\mbmss{x} \quad\, A_{i_1}^{n_1} \atop f(\mbmss{x})\,f(A_{i_1}^{n_1})} {\cdots\,\,\,\, A_{i_r}^{n_r}\,\atop \,\,\,\cdots \, f(A_{i_r}^{n_r})}\big\rangle$ can be omitted as follows.
 We also write 
\[ \S_{\langle {\mbmty{x}\atop \mbmty{\xi}}{A_{i_1}^{n_1} \atop \alpha_1}{\cdots\atop \cdots} {A_{i_r}^{n_r}\atop \alpha_r}\rangle}(H(\mbm{x},A_{i_1}^{n_1} ,\ldots,A_{i_r}^{n_r}))\] for $\S_{f}(H(\mbm{x},A_{i_1}^{n_1} ,\ldots,A_{i_r}^{n_r}))$ if $\mbm{\xi}=f(\mbm{x})$, $\alpha_1=f(A_{i_r}^{n_1})$, \ldots, $\alpha_r=f(A_{i_r}^{n_r})$, and the like.
\end{agree}

\vspace{0.2cm}\noindent {\bf The assignments in ${\rm assgn}_H(\S)$.}
Let $\S= (J_n)_{n\geq 0}$ be in ${\sf struc}_{\rm pred}^{({\rm m})}(J_0)$, $H$ be in ${\cal L}^{(2)}$, and
 $x_1,\ldots,x_s$ and $A_{i_1}^{n_1}, \ldots,A_{i_r}^{n_r}$ be the only variables that may occur free in $H$. Then, let ${\rm assgn}_H(\S)$ be the set of all {\em restricted {\em $\mbm{\sigma}^{(2)}$-assignments $(f_n)_{n\leq \max\{n_1,\ldots,n_r\}}$ in $\S$} with $f_0:\{x_i\mid i\geq 0\}\to J_0$ and $f_n:\{A^n_i\mid i\geq 0\}\to J_n$ for $n\geq 1$ with $n\leq \max\{n_1,\ldots,n_r\}$}.

\begin{defi}[Valid formulas, models, $\models_f$, and $\models$] 
Let $f$ be an arbitrary assignment in ${\rm assgn}_H(\S)$ and in ${\rm assgn}(\S)$, respectively.
We say that the formula {\em $ H$ is true in $\S$ under $f$} and {\em has the truth value $true$ in $\S$ under $f$} and we write $\S_f\models H$ or $\S\models_fH$ in analogy to the notation in \cite[Definition\,\,2]{Vaeae} if $\S_f(H)=true$. 
$H$ is {\em valid in $\S$} iff we have $\S\models_fH$ for each assignment $f$ in ${\rm assgn}_H(\S)$. We say that
 $\S$ is a {\em model of $H$} and write $\S \models H$ if $H$ is valid in $\S$. We write $\models H$ if $\S \models H$ holds for all $\S$ in ${\sf struc}_{\rm pred}^{({\rm m})}$. For any subset ${\cal G}\subseteq {\cal L}^{(2)}$, we say that $\S$ is a {\em model of ${\cal G}$} and we  write $\S\models {\cal G}$ if each $H\in {\cal G}$ is valid in $\S$. ${\cal G}\models H$ means that  $\S\models H$ holds for all predicate structures $\S$ with $\S\models {\cal G}$.
 \end{defi}

\subsection{Predicate structures and definability in PL\,II}\label{Section_PL_II}

Let us  start with  the axioms in $axa^{(2)}$, $ axi^{(2)}$, and $ext^{(2)}$, and let us then add further axioms, i.e.,  the axioms in $komp^{(2)}$.  For details see, e.g., \cite{1a,1b,1c}.
For the axioms of propositional logic see also \cite{HiBe34} and in particular Frege's Begriffsschrift \cite{Frege}.

\vspace{0.2cm} \noindent {\bf The set $ax^{(2)}$ of axioms of PL\,II.} 
 Let $ax^{(2)}$ be the {\em axiomatic system of {\rm PL\,II}} consisting of all axioms in $axa^{(2)}\cup axi^{(2)}\cup ext^{(2)}$.

\begin{list}{(A\arabic{li})}
{\usecounter{li}} 
\item 
$axa^{(2)}$:
The axioms of propositional logic for PL\,II 
(cf.\,\cite{1a})
\item 
$axi^{(2)}$:\,
The axioms of identity theory for PL\,II 
(cf.\,\cite{1b} and \cite{1c})
\item 
$ext^{(2)}$:\hspace{0.08cm}
The axioms of extensionality for PL\,II (cf.\,\cite{1c})
\end{list}

\begin{theorem}[Asser \normalfont{\cite{1c}}]$\!\!\!$ The predicate structures are --- up to
isomor\-phisms --- the only structures of signature $\mbm{\sigma}^{(2)}$ that are models 
of $ax^{(2)}$. 
\end{theorem}

\begin{defi}[Second-order definability] Let $\S$ be in ${\sf struc}_{\rm pred}^{({\rm m})}(J_0)$, $H$ be a second-order formula in which the variables $x_1,\ldots,x_n$ may occur only free and $f$ be an assignment in ${\rm assgn}(\S)$ or in ${\rm assgn}_H(\S)$. A predicate $\alpha: J_0^n\to \{true, false\}$ is {\em second-order definable over $\S$ with the help of $H$ and $f$} if we have $\alpha(\mbm{\xi})= \S_{f\langle {\mbmty{x}\atop\mbmty{\xi}}\rangle}(H(x_1,\ldots,x_n))$ for all $ \mbm{\xi}\in J_0^n$. 
\end{defi}

Let $x_1,\ldots,x_n$, $x_{j_1},\ldots,x_{j_s}$, and $A_{i_1}^{n_1},\ldots,A_{i_r}^{n_r}$ be the only variables that may occur free in $H$. Then, we can also say that the $n$-ary predicate $\alpha$ is {\em second-order definable over $\S$ with the help of  $f( x_{j_1} ),\ldots,f(x_{j_s})$ and  $f(A_{i_1}^{n_1}), \ldots,f(A_{i_r}^{n_r})$} if we have $\alpha(\mbm{\xi})= \S_{f\langle {\mbmty{x}\atop\mbmty{\xi}}\rangle}(H(x_1,\ldots,x_n,x_{j_1},\ldots,x_{j_s},A_{i_1}^{n_i} ,\ldots,A_{i_r}^{n_r}))$ for all $ \mbm{\xi}\in J_0^n$. 

For $H\in {\cal L}^{(2)}_{x_{i_1},\ldots,x_{i_n}}$, let 
\[\alpha_{\S,H,x_{i_1},\ldots,x_{i_n}, f}(\xi_{i_1},\ldots,\xi_{i_n}) =\S[H(x_{i_1},\ldots,x_{i_n}), f\langle {{x_{i_1}\atop \xi_{i_1}}{\cdots\atop\cdots} {x_{i_n}\atop \xi_{i_n}}\rangle}].\]

In the following, we will also use set-theoretic operations to define predicates.

\begin{agree}\label{AgreeNamesForPred} For each predicate $\alpha$, we also use the notation $\alpha$ for the predicate relation $\widetilde \alpha=\{\mbm{\xi}\in I^n\mid \alpha(\mbm{\xi})=true\}$ whenever a confusion is not possible and the context allows to recognize whether the used notation $\alpha$ denotes the relation $\widetilde \alpha$ or the predicate $\alpha$ itself. 
\end{agree}

Moreover, for $n$-ary predicates $\alpha$, the abbreviation $\alpha \subseteq I^n$ stands for $\widetilde\alpha \subseteq I^n$ and $\alpha = I^n$ stands for the statement $\widetilde\alpha = I^n$. We say a predicate $\alpha$ is {\em finite} if $\widetilde\alpha$ is finite, and so on. 

\subsection{Henkin structures and closeness with respect to definability in HPL}\label{Section_HPL}

Let $\S$ be any structure $(J_n)_{n\geq 0}$ in ${\sf struc}_{\rm pred}^{({\rm m})}(J_0)$. 
Let $Att_\S^n$ be the operator that provides $\alpha_{\S,H,\mbmss{x}, f}$ for any $H$ in ${\cal L}^{(2)}_{\mbmss{x}}$ and any assignment $f$ in ${\rm assgn}(\S)$ or ${\rm assgn}_H(\S)$. Shortly we write $Att_\S^n[H, f]=\alpha_{\S,H,\mbmss{x}, f}$ for expressing this fact. 
For any $n$, any $H(\mbm{ x})$, and any $f$, there are two possibilities. Either $Att_\S^n[H, f]$ belongs to $J_n$ or it does not belong to $J_n$. 

For $H\in {\cal L}^{(2)}_{x_{i_1},\ldots,x_{i_n}}$, let $Att_\S[H;x_{i_1},\ldots,x_{i_n} ;f]=\alpha_{\S,H,x_{i_1},\ldots,x_{i_n}, f}$. Thus, for any $H$ in ${\cal L}^{(2)}_{\mbmss{x}}$, we have \[Att_\S^n[H, f]=Att_\S[H;x_1,\ldots,x_n ;f].\]

\begin{defi}[Closeness with respect to definability] 
Let a predicate structure $\S$ be {\em closed under $Att_\S^n$} if the predicate $Att_\S^n[H, f]$ belongs to $\S$ for all formulas $H$ and all $f$ in ${\rm assgn}_H(\S)$. 
We say that {\em $\S$ is closed with respect to second-order definability} if it is closed under $Att_\S^n$ for all $n\geq 1$.
\end{defi}

When $\S$ is closed under $Att_\S^n$ and $A_{0}^{n}$ does not occur in $H(\mbm{x})$, then $\S\models _f\exists A_{0}^{n} \forall \mbm{x} (A_{0}^{n}\mbm{x}\leftrightarrow H(\mbm{x}))$ holds and the value $\S_f(A_{0}^{n}\mbm{x}\leftrightarrow H(\mbm{x}))$ is $true$ for $f(A_{0}^{n})=Att_\S^n[H, f]$. Consequently, the predicate structures that are closed with respect to second-order definability are models of the axioms of comprehension.
In \cite{1c}, the predicate structures of this special type are introduced and investigated. Thus, we also call them  {\em Henkin-Asser structures}. 

\begin{defi}[The Henkin structures of second order] \hfill A predicate \linebreak structure that is also a model of all axioms resulting from the schema given by (A4) is said to be a {\em Henkin structure} ({\em of second order}) or a {\em Henkin-Asser structure}.  \begin{list}{(A\arabic{li})}
{\usecounter{li}\setcounter{li}{3}}
\item $komp^{(2)} $: 
The axioms of comprehension for {\rm HPL} (cf.\,\cite{1c})
\end{list}\end{defi}

Let $^{h}ax^{(2)}$ be the {\em axiomatic system of {\rm HPL}} consisting of all axioms in $ax^{(2)}\cup komp^{(2)} $ where $komp^{(2)}$ is the set of all axioms given by (A4).

\begin{theorem}[Asser \normalfont{\cite{1c}}] The Henkin structures of second order are --- up to isomor\-phisms --- the only structures of signature $\mbm{\sigma}^{(2)}$ that are models of $^{h}ax^{(2)}$. 
\end{theorem}
Thus,  we have the following corollary. 
\begin{corollary}\label{ClosenessDefinability}
Each predicate structure $\S$ that is closed with respect to the second-order definability is a Henkin-Asser structure. 
\end{corollary}
We will deal with models of  PL II and HPL.

\vspace{0.3cm}

\fbox{\parbox{11.3cm}{
\sf \small
Second-Order Models of PL II
\begin{itemize}
\item Models of $axa^{(2)}\cup axi^{(2)}\cup ext^{(2)} $ (up to isomorphisms): 

\quad predicate structures
\end{itemize}

\vspace{0.1cm}
Second-Order Models of HPL 
\begin{itemize}
\item Models of $^{h}ax^{(2)}$ (up to isomorphisms): 

\quad Henkin structures 
\end{itemize}
}}

\vspace{0.5cm} 
By analogy with Henkin's Theorem of Completeness in \cite[Theorem\,2]{Henk}, Asser formulated the Hauptsatz für das $\mathfrak{H}^{(2)}$-Folgern in the restricted framework as follows (where $\mathfrak{H}^{(2)}$-Folgern is the process of generating semantic consequences by applying $\models_h$ defined by ${\cal G}\models_h H$ iff  $^hax^{(2)}\cup {\cal G}\models H$). For details see \cite{1c}. 

\begin{theorem}[The completeness theorem{\normalfont,} Asser {\normalfont \cite{1c}}]\label{Completeness}For every  second-order formula $H$ in ${\cal L} ^{(2)}$, we have $^{h}ax^{(2)}\models H$ if and only if $^{h}ax^{(2)}\vdash H$. 
This means that $H$ is generally valid in {\rm HPL} if and only if $H$ is derivable from the axioms given by (A1), (A2), (A3), and (A4) with respect to the usual (syntactical) rules of inference for $\vdash$.
\end{theorem}

\subsection{Permutations of individuals and their extensions} 
For any non-empty set $I$, each {\em permutation of $I $} is a bijection from $I$ to $I$. Let $({\rm perm}(I),\circ)$ and $(\G_{\sf 1}^I,\circ)$ (shortly denoted by ${\rm perm}(I)$ or $\G_{\sf 1}^I$) be the group of all permutations of $I$ with the usual function composition $\circ$ as binary group operation. Let $I$ be any arbitrary non-empty set of individuals and $\pi$ be a permutation in $\G_{\sf 1}^I$. Then, for all $\mbm{\xi} \in I^n$, we will use the abbreviation $\pi(\mbm{\xi})$ defined by $\pi(\mbm{\xi})= (\pi(\xi_{1}), \ldots,\pi(\xi_{n}))$. Moreover, for each $n$-ary predicate $\alpha\in {\rm pred}_n(I)$, let $\alpha^\pi$ be the $n$-ary predicate satisfying \[\alpha^\pi(\mbm{\xi})=\alpha(\pi^{-1}(\mbm{\xi}))\] for all $\mbm{\xi}\in I^n$.\footnote{The definition of $\alpha^\pi$ given in \cite{1c} was modified.} Thus, we have $\widetilde{\alpha^{\pi}}= \{\pi(\mbm{\xi})\mid \mbm{\xi}\in \widetilde\alpha\}$. Let $\G$ be any subgroup of $\G_{\sf 1} ^I$. The group $(\G,\circ)$ (shortly, $\G$) can be any group of bijections from $I$ to $I$ or a group of automorphisms of a given basic structure ${\sf S}$ with $I$ as universe. 
It is easy to extend every permutation $\pi$ in $\G$ uniquely to a permutation $\pi^*$ of $I\cup {\rm pred}(I)$ by defining
\[\pi^*(\alpha)=\alpha^\pi \quad \mbox{ for all }\alpha \in {\rm pred}(I) \quad \mbox{ and }\quad \pi^*(\xi)=\pi(\xi) \quad \mbox{ for all }\xi\in I. \] Let $\G^*= \{ \pi^* \mid \pi\in \G\}$.

Now, we define a second extension of the permutations by expanding the domain to assignment functions. Let $I$ be any non-empty set, $\S$ be any predicate structure $(J_n)_{n\geq 0}$ in ${\sf struc}^{({\rm m})}_{\rm pred}( I)$ which means $I=J_0$, and $\pi\in \G_{\sf 1}^I$. Then, it shall be possible to apply a function from ${\rm assgn}(\S)$ to ${\rm assgn}(\S_{\sf 1}^{I})$ determined by $\pi$. 
For every $f$ in ${\rm assgn}(\S)$, let $f^{\pi}$ be the assignment in ${\rm assgn}(\S_{\sf 1}^{I})$ defined as follows.\footnote{The definition of $f^{\pi}$ given in \cite{1c} was modified.} 

\vspace{0.3cm}
\begin{tabular}{ll}
 $f_0^{\pi}(x_i)$ & $=_{\rm df}\pi(f_0(x_i)) $\hfill{ \hspace{2.8cm} (for all $f_0:\{x_i\mid i\geq 0\}\to I$)}\\[1.5ex]
 $f_n^{\pi}(A_i^n)$ & $=_{\rm df}(f_n(A_i^n))^{\pi} $\hfill{ (for all $f_n:\{A^n_i\mid i\geq 0\}\to J_n$, $n\geq 1$)}\\[1.5ex]
 $f^{\pi}$ & $=_{\rm df}(f_n^{\pi})_{n\geq 0}$\hfill{ (for all $f=(f_n)_{n\geq 0}$)}
\end{tabular}
\vspace{0.3cm}

\noindent For $f_0:\{x_i\mid i\geq 0\}\to I$ we get $f_0^\pi:\{x_i\mid i\geq 0\}\to I$. For all $f_n:\{A^n_i\mid i\geq 0\}\to J_n$, we get $f_n^\pi:\{A^n_i\mid i\geq 0\}\to {\rm pred}_n(I)$.

\begin{lemma}[Closeness of $ {\rm assgn}(\S)$ with respect to $\pi$] Let $\S=(J_{n})_{n\geq 0}$ and $\pi\in \G_{\sf 1}^{J_0} $. If all images (all function values) of $ \pi^ *$ are in $\S$, then there holds $f^{\pi}\in {\rm assgn}(\S)$.
\end{lemma}
 Such an assignment $f^{\pi}$ has the following three properties. 
\begin{lemma}[Asser \normalfont{\cite{1c}}]\label{HSAsser} Let $\S$ be the predicate structure $(J_n)_{n\geq 0}$ and let $\pi$ be a permutation in $\G_{\sf 1}^{J_0} $ such that all images of $ \pi^ *$ are in $\S$. Let $f$ be an assignment in $\S$ and $H$ be any formula in ${\cal L}^{(2)}$. 
 Then, we have {\em (1)} and we have the properties {\em (2)} and {\em (3)} if $H$ is in ${\cal L}^{(2)}_{x_{1},\ldots,x_{n}}$.
\begin{list}{\normalfont (\arabic{li})}{\usecounter{li}\labelwidth0.5cm \leftmargin0cm \itemsep2pt plus1pt\topsep1pt plus1pt minus1pt\labelsep4pt \parsep0.5pt plus0.1pt minus0.1pt \itemindent0.8cm}
\item\label{HSAsser1} $\S_{f^{\pi}}(H)=\S_f(H)$.
\item\label{HSAsser2} $\alpha_{\S,H,\mbmss{x}, f^{\pi}}=\alpha_{\S,H,\mbmss{x}, f} ^{\pi}$.
\item\label{HSAsser3} $\alpha_{\S,H,\mbmss{x}, f^{\pi}}$ is definable over $\S$ if $\alpha_{\S,H,\mbmss{x}, f} $ is definable over $\S$. 
\end{list}
\end{lemma} 
\vspace{0.2cm}

By Lemma \ref{HSAsser} (\ref{HSAsser3}), each extension $\pi^*$ of a permutation $\pi$ in a subgroup $\G$ of $perm(I)$ transforms the predicate $\alpha_{\S,H,\mbmss{x}, f}$ definable by $H$, individuals $\xi_{n+1},\ldots,\xi_{n+s}$, and predicates $\alpha_{i_1}, \ldots,\alpha_{i_r}$ into a predicate definable by $H$. If the mentioned individuals and predicates are mapped to themselves by means of $\pi$, then the relation $\widetilde \alpha_{\S,H,\mbmss{x}, f}$ is also transformed into itself by $\pi$. Since we are especially interested in constructing structures satisfying the axioms given by (A4), we consider so-called symmetry subgroups of a group $\G$ of permutations, stabilizer subgroups, and supports for predicates and discuss their definability by means of second-order formulas.

\subsection{The basic Fraenkel model of second order}\label{AbschnFrae_mod}

For constructing Henkin structures, we can follow \cite{1c} where the fundamental papers \cite{M38}, \cite{M39}, and \cite{FA22a} are cited in this context. 

 {\bf Symmetry subgroups and invariance}\label{SectionSymmetrySubgroups}
Let $I$ be any non-empty set and $\G$ be a subgroup of $perm(I)$. For each predicate $\alpha \in {\rm pred}(I)$, let ${\rm sym}_{\G}(\alpha)$ be defined by
\[{\rm sym}_{\G}(\alpha)=\{ \pi \in {\G}\mid \alpha^{\pi}=\alpha\}\]
 and let ${\rm sym}_{\G}^*(\alpha)=_{\rm df} ({\rm sym}_{\G}(\alpha))^*$. Then, ${\rm sym}_{\G}(\alpha)$ is a group called the {\em symmetry subgroup of $\G$ with respect to $\alpha$}.  Whereas $\G^*$ is a subgroup of $(perm(I))^*$, ${\rm sym}_{\G}^*(\alpha)$
is a subgroup of $\G^*$ that also acts on $I\cup {\rm pred}(I)$ from the left.  ${\rm sym}_{\G}^*(\alpha)$ is the stabilizer subgroup (i.\,e.\,\,the isotropy subgroup) of $\G^*$ with respect to $\alpha$. Moreover, let $perm^*(I)=(perm(I))^*$. Thus, for all $\alpha \in {\rm pred}(I)$, we have $\pi^*(\alpha)= \alpha ^\pi$ for $\pi^*\in perm^*(I)$ and $\pi\in perm(I)$, respectively.

For any predicate $\alpha\in {\rm pred}(I)$, ${\rm sym}_{\G}(\alpha)$ is the symmetry subgroup of all permutations $\pi$ in $\G$ whose extensions $\pi^*$ transform $\alpha$ into itself. A predicate $\alpha$ is {\em fixed with respect to a permutation $\pi^* \in\G^*$} if $\pi^*(\alpha)=\alpha$. Thus, ${\rm sym}_{\G}^*(\alpha)$ is the subgroup of all permutations $\pi^* \in {\G}^*$ with respect to which $\alpha$ is fixed. Note that this means that $\alpha$ is a fixed point of all permutations in ${\rm sym}_{\G}^*(\alpha)$. Moreover, an $n$-ary relation $\mbm{r}\subseteq I^n$ is {\em invariant with respect to a permutation $\pi \in \G$} if $\mbm{\xi} \in \mbm{r}$ holds if and only if $\pi( \mbm{\xi}) \in \mbm{r}$ holds. Thus, ${\rm sym}_{\G}(\alpha)$ is also the subgroup of all permutations $\pi \in {\G}$ with respect to which $\widetilde\alpha$ is invariant. This means that $\widetilde \alpha$ is invariant under ${\rm sym}_{\G}(\alpha)$ and that ${\rm sym}_{\G}(\alpha)$ is the subgroup of $\G$ consisting of all automorphisms of the structure $(I;\emptyset;\emptyset;\widetilde \alpha)$ (that contains no constants and functions). 

\begin{agree}
In the following, we also omit the asterisk after $\pi$ in $\pi^*$ and write then also $\pi$ if we consider $\pi^*$.
\end{agree} 

\noindent{\bf Finite supports for predicates.} 
For any subgroup $\G$ of $perm(I)$ and any finite $P\subseteq I$, let ${\G}( P)$ be the subgroup containing all $\pi \in {\G}$ which map each element in $P$ to itself. Another notation for such a group could --- in accordance with the notations used in \cite{Jech73} --- be ${\rm fix}_{\G}(P)$. More precisely, let
\[{\G}( P)=\{\pi \in \G\mid \pi(x)=x \mbox{ for all } x \in P\} \mbox{ and } {\rm fix}_{\G}(P)={\G}( P).\]

The basic Fraenkel model of second order, in the following denoted by $\S_0$, can be defined by using an infinite set $\bbbn$ of individuals, the group $perm(I)$, and finite supports.

\begin{agree}[Our set $\bbbn$]Here, let $\bbbn$ be an arbitrary infinite set of individuals. 
\end{agree} 

\begin{agree}[The model $\S(I,\G,\F)$] Here, we only use second-order permutation models of HPL denoted by $\S(I,\G,\F)$ where $I$ is the non-empty individual domain, $\G$ is a subgroup of $perm(I)$, and $\F$ is the following normal filter $\F_{\sf \! 0}(I,\G)$ on $\G$.  $\F_{\sf \! 1}(I, \G)$ contains all subgroups of $\G$. $\I_0^I$ is a normal ideal containing all finite subsets of $I$.
\end{agree} 
We will use the definitions
\[\F_{\sf \! 0}(I, \G)=_{\rm df}\{\H \in \F_{\sf \! 1}(I, \G)\mid (\exists P\subseteq I)( P \mbox{ \rm  is a finite set} \,\,\,\&\,\, \H\supseteq \G(P))\}\]
and
\[\F_{\sf \! 0}(I, \G)=_{\rm df}\{\H \in \F_{\sf \! 1}(I, \G)\mid (\exists P\in \I_0^I)( \H\supseteq \G(P))\}.\]

We say that the predicate $\alpha$ is {\em symmetric with respect to $\F$} if the symmetry subgroup ${\rm sym}_{\G}(\alpha)$ is in $ \F$.
For any $n\geq1$, let $J_{n}(I,\G,\F)=\{\alpha \in {\rm pred}_n (I)\mid {\rm sym}_{\G}(\alpha)\in \F\}$ and let $\S(I,\G,\F)$ be the predicate structure $ (J_{n})_{n\geq 0}$ with the domains $J_{0}=I$ and $J_{n}=J_{n}(I,\G,\F)$ for $n\geq 1$. Thus, for every $n\geq 1$, $J_{n}(I,{\G,\F})$ is the set of all $n$-ary predicates that are symmetric with respect to $\F$. The following closeness property shows the usefulness of the normal filter $\F_{\sf \! 0}(I, \G)$ (cf.\,Asser \cite{1c}).

\begin{proposition}\label{Closen_Filt} The structure $\S(I,\G,\F)$ is closed with respect to all permutations $\pi \in \G$.
\end{proposition}

Let $\S_0$ be the basic Fraenkel model of second order $ \S(\bbbn,\G_1^\bbbn,\F_{\sf \! 0}(\bbbn, \G_1^\bbbn))$.

\vspace{0.2cm} \noindent Now, for  any subgroups $\G$ of $\G_1^I$ and any normal filter $\F$ on $\G$ we will characterize the stabilizers for the truth values of formulas by the following generalization of Hilfssatz 4 in {\normalfont \cite{1c}}. In particular, we can say that the subgroup $ {\rm sym}^*_{\G}(\alpha)$ is {\em the stabilizer subgroup of $\G^*$ with respect to a formula $H$ under an assignment $f$ in $\S(I,{\G,\F})$} if $\alpha=\alpha_{\S,H,\mbmss{x}, f}$. Proposition \ref{Closen_Filt}, together with Lemma \ref{HSAsser}\,(\ref{HSAsser2}), implies the following proposition that can be proved by using the observations discussed for finite supports in \cite{1c}. 

\setcounter{equation}{0}
\begin{proposition}[Stabilizers for formulas]\label{AsserHS4} Let $\S=\S(I,{\G,\F})$ and $f\in {\rm assgn}(\S)$. Moreover, let $H(\mbm{x})$ be in ${\cal L}^{(2)}_{\mbmss{x}}$ and $r,s\geq 0$. If $s>0$, then let $x_{i_1},\ldots,x_{i_s}$ ($ i_1,\ldots,i_s>n$) be the individual variables that may additionally occur free in $H(\mbm{x})$. If $r>0$, then let $A_{j_1}^{n_1}, \ldots, A_{j_r}^{n_r}$ be the predicate variables that may occur free in $H(\mbm{x})$. Let $H(\mbm{x})$ does not contain further free variables. Let $\H_{r+1}= \G(\{f(x_{i_1}) , \ldots, f(x_{i_s})\})$ if $s>0$ and $\H_{r+1}= \G$ if $s=0$. If $r>0$, then, for every $k\in\{1,\ldots,r\}$, let $A_k$ stand for $A_{j_k}^{n_k}$ and let $\H_k $ be a subgroup in $\F$ such that ${\rm sym}_{\G}(\alpha_k)\supseteq \H_k$ holds for $\alpha_k=f(A_k)$. For each $l\in \{1,\ldots,r+1\}$, let $\G_l$ be the subgroup in $\F$ given by
 $\G_l= \H_{l}\cap \cdots \cap\H_{r+1}$. Then, we have 
\begin{equation}\label{Stab_1}\S_f(H(\mbm{x},A_1,\ldots, A_{m}))=\S_{f\langle {\mbmty{x}\atop (\pi(f(\mbmty{x}))} {A_1\atop (\alpha_1)^ \pi} {\cdots\atop\cdots} {A_{m}\atop (\alpha_{m})^ \pi} \rangle}(H(\mbm{x},A_1,\ldots, A_{m}))\end{equation} for each each $m\in \{1,\ldots, r\}$ and $\pi \in \G_{m+1}$ if $r>0$ and
 \begin{equation}\label{Stab_2}{\rm sym}_{\G}(\alpha_{\S,H,\mbmss{x}, f})\supseteq \G_1\end{equation} in any case where $r\geq 0$.
\end{proposition}

 \noindent Equation (\ref{Stab_2}) in Proposition \ref{AsserHS4} means that, for each predicate definable in a predicate structure $\S(I,{\G,\F})$, the symmetry subgroup of $\G$ with respect to this predicate belongs to the considered filter $\F$. This guarantees the closeness of $\S(I,{\G,\F})$ with respect to the definability by means of second-order formulas. This means that $\S(I,{\G,\F})$ is closed under $Att_\S^n$ for all $n\geq 1$.  A consequence is that  structures constructed such as $\S_0$ are  models of the axioms given in (A4).

\begin{corollary}[Asser {\normalfont \cite{1c}}]\label{AsserHenkin} The predicate structure $\S(I,\G,\F)$ is a  Henkin structure.
\end{corollary}

We want to continue our investigations. We  consider  $\S_0$ in order to prove the independence of the well-ordering theorem from a more general principle of choice in second-order logic HPL.

\section{Classical second-order versions of AC}\label{AbschnittAC}\label{SectionTheAxiom} 
\setcounter{satz}{0}

In formalizing several principles of choice and the principle of well-ordering we will also use new variables for predicates whose arities are determined
by the types of the corresponding tuples of variables in the formulas. For given $n\geq 1$ and $m\geq 1$, let $A$ and $B$ be variables of sort $n$, let $C, C_1, C_2$, and $D$ be of sort $m$, let $R$ and $S$ be of sort $n+m$, and let $T$ be of sort $2n$. 

We start with the definition of the principles of choice of second-order logic in a form discussed in \cite{1c}. First, we define Ackermann's principle (cf.\,\,\cite{Ack111}) as introduced by means of Church's $\lambda$-notation (cf.\,\cite{Church32}) in \cite{1c}. For this reason, let $n\geq 1$, $m\geq 1$, and $H(\mbm{x},D) $ be any second-order formula in ${\cal L}^{(2)}_{\mbmss{x},D}$ in which the variables $x_1, \ldots, x_n$, and $ D$ occur only free. Moreover, for any predicate structure $\S= (J_n)_{n\geq 0}$, any $(n+m)$-ary predicate $\sigma\in J_{n+m}$, and any $\mbm{\xi}\in J_0^n$, let $\S_{\langle {\mbmty{x}\atop \mbmty{\xi}}{ S\atop\sigma}\rangle} (\lambda\mbm{y}.S\mbm{x}\mbm{y})$ be the $m$-ary projection $\sigma_{\mbmss{\xi}}\in J_{m}$ satisfying $\sigma_{\mbmss{\xi}}(\mbm{\eta}) =\sigma(\mbm{\xi}\,.\,\mbm{\eta}) $ for every $\mbm{\eta}\in J_0^m$.
Then, roughly speaking, an important principle of choice states that the condition $\S\models_f \forall \mbm{x} \exists D H(\mbm{x},D)$ implies the existence of an $(n+m)$-ary predicate $\sigma\in J_{n+m}$ such that, for each $n$-tuple $\mbm{\xi}\in J_0^n$, $\S\models_{f\langle{\mbmty{x}\atop \mbmty{\xi}}{D\atop{\sigma}_{\mbmty{\xi}}}\rangle} H(\mbm{x},D)$ holds.

\vspace{0.3cm}
\noindent{\bf The Ackermann axioms $choice^{n,m}(H)$}
\nopagebreak 

\noindent\fbox{\parbox{11.7cm}{

\vspace*{0.1cm}
$choice^{n,m}(H)=_{\rm df}  \forall \mbm{x} \exists D\, H(\mbm{x},D) \to \exists S \forall \mbm{x}\, H(\mbm{x}, \lambda\mbm{y}.S\mbm{x}\mbm{y})
$
}} 

\vspace*{0.3cm}
\noindent
 Let $choice^{(2)}$ be the set of all Ackermann axioms.

\vspace{0.3cm}
\noindent{\bf The set $choice^{(2)}$ of all Ackermann axioms in PL\,II} 

\nopagebreak 

\noindent\fbox{\parbox{11.7cm}{

\vspace*{0.1cm}

$\begin{array}{ll}choice^{(2)}&=_{\rm df}\bigcup_{n,m\geq 1}choice^{n,m}\\
choice^{n,m}&=_{\rm df}\{choice^{n,m}(H)\mid H \mbox{ is a formula in ${\cal L}^{(2)}_{\mbmss{x},D}$} \}\\
\end{array}$
}} 

\vspace*{0.3cm}

According to the axiom schema given by (A4), there is a further possibility to formulate these principles of choice without $\lambda$-terms. 

\vspace{0.3cm}
\noindent{\bf The Ackermann axioms $choice^{n,m}_h(H)$} 

\nopagebreak 

\noindent\fbox{\parbox{11.7cm}{

\vspace*{0.1cm}

$ choice_h^{n,m}(H)=_{\rm df} 
 \forall \mbm{x} \exists D H(\mbm{x},D) \to \exists S \forall \mbm{x} \exists D (\forall\mbm{y}(D \mbm{y} \leftrightarrow S\mbm{x}\mbm{y}) \land H(\mbm{x}, D )) 
$
}} 

\vspace*{0.3cm}
\noindent Let $choice_h^{(2)}$ be the set of all these Ackermann axioms.

\vspace{0.3cm}
\noindent{\bf The set $choice_h^{(2)}$ of all Ackermann axioms in HPL} 
\nopagebreak 

\noindent\fbox{\parbox{11.7cm}{

\vspace*{0.1cm}
$\begin{array}{ll}choice_h^{(2)}&=_{\rm df}\bigcup_{n,m\geq 1}choice_h^{n,m}\\
choice_h^{n,m}&=_{\rm df}\{choice_h^{n,m}(H)\mid H \mbox{ is a formula in ${\cal L}^{(2)}_{\mbmss{x},D}$} \}\\
\end{array}$
}} 

\vspace*{0.3cm}
\noindent For any $H(\mbm{x},D)$ in ${\cal L}^{(2)}_{\mbmss{x},D}$, $ choice^{n,m} (H)$ is {\rm HPL}-equivalent to $ choice_h^{n,m} (H) $.

If all variables in the lists $\mbm{x}$ and $\mbm{y}$ occur only free in the second-order formula $ H$, then the following formula is a generalization of formulas such as the axioms that are given for $m=1$ in \cite{HiAc59}.

\vspace{0.3cm}
\noindent{\bf The Hilbert-Ackermann axioms $AC^{n,m}(H)$ of second order}

\nopagebreak 

\noindent\fbox{\parbox{11.8cm}{

\vspace*{0.1cm}
$AC^{n,m}(H)\!=_{\rm df}\!\forall A \exists S
(\forall \mbm{x}
 ( A\mbm{x}\! \leftrightarrow \! \exists \mbm{y}
 H(\mbm{x},\!\mbm{y})) \!\to \forall \mbm{x}
 (A\mbm{x}\!\to\! \exists!! \mbm{y}
 (H(\mbm{x},\!\mbm{y})
 \land S\mbm{x}\mbm{y})))$
}} 
\vspace{0.3cm}

By replacing $H(\mbm{x},\mbm{y})$ in $AC^{n,m}(H)$ by $R\mbm{x}\mbm{y}$, we get further formulations of {\em Zermelo's axioms of choice} in second-order logic. $AC^{n,m}(R\mbm{x}\mbm{y})$ was denoted by $Ch_ {n,m}$ in \cite{Gass84} and will be denoted by $AC^{n,m}$ in the following.

\vspace{0.3cm}
\noindent{\bf The Zermelo-Asser axioms $AC^{n,m}$ of second order} 

\nopagebreak 

\noindent\fbox{\parbox{11.8cm}{

\vspace{0.1cm}
$AC^{n,m}=_{\rm df}\forall A \forall R \exists S
(\forall \mbm{x}
 ( A\mbm{x} \leftrightarrow \exists \mbm{y}
 R\mbm{x}\mbm{y}) \to \forall \mbm{x}
 (A\mbm{x}\to \exists!! \mbm{y}
 (R\mbm{x}\mbm{y}
 \land S\mbm{x}\mbm{y})))$
}}
\vspace{0.2cm}

\noindent This form is more general than the formula $Ch_ {n,m}^{(2)}$ introduced by Asser  \cite{1c}. $R$ can now stand for an arbitrary $(n+m)$-ary relation and the main statement is restricted to the {\em $n$-ary domain} $A$ of $R$. In \cite{1c}, $AC^{1,1}$  was introduced as a further variant of the axiom of choice and $AC^{n,m}$ was called the {\em $(m,n)$-variant of $AC^{1,1}$}. $AC^{n,m}$ includes the relativization of $Ch_ {n,m}^{(2)}$ to each arbitrary $n$-ary predicate for which the variable $A$ can stand. In {\rm HPL}, it is equivalent to Asser's formula $Ch_ {n,m}^{(2)}$. 
It is clear that the Axiom of Choice --- given by Jech \cite{Jech73} in the well-known formulation --- would hold in (our metatheory) ZF if $AC^{1,1}$ would hold in any standard structure of second order. 

The following formulas --- given in \cite{1c} and denoted by $AC_{n,m}^{\rm Disj}$ in \cite{Gass84} and by $AC_*^{n,m}$ in \cite{Gass94} --- are a result of the formalization of {\em Russell's axioms of choice} in second-order logic.

\vspace{0.3cm}
\noindent{\bf The Russell-Asser axioms $AC_*^{n,m}$ of second order}

\nopagebreak

\noindent\fbox{\parbox{11.8cm}{

\vspace{0.1cm} 
\noindent $AC_*^{n,m}=_{\rm df}\forall A \forall R \exists S
(\forall \mbm{x}
 ( A\mbm{x} \! \leftrightarrow \! \exists \mbm{y}
 R\mbm{x}\mbm{y} ) $

\hfill $
 \land \forall \mbm{x}_{1} \forall \mbm{x}_{2}
 ( A\mbm{x}_1 \land A\mbm{x}_2 \land \mbm{x}_{1}\! \neq \!\mbm{x}_{2}
 \to \lnot \exists \mbm{y}
 (R\mbm{x}_{1}\mbm{y}
 \land R\mbm{x}_{2}\mbm{y}))$
 
 \hspace{3cm}$
 \to \forall \mbm{x}(A\mbm{x}
 \to \exists !!\mbm{y}(R\mbm{x}\mbm{y}
 \land S\mbm{x}\mbm{y})))$}}

\vspace{0.2cm} 

 In \cite{1c}, Asser also gave a further formulation. After a small modification, we get the axioms for any $H$ in ${\cal L}_C^{(2)}$. 

\vspace{0.3cm}
\noindent{\bf The Asser axioms $choice_*^{m}(H)$} 

\nopagebreak

\noindent\fbox{\parbox{11.8cm}{

\vspace{0.1cm} 
$choice_*^{m}(H)=_{\rm df}\forall C(H(C)\to\exists \mbm{y}C\mbm{y}) $

\hfill $\land\,\, \forall C_1\forall C_2 (H(C_1)\land H(C_2)\land C_1\neq C_2\,
\,\to\, \neg \exists \mbm{y}(C_1\mbm{y}\land C_2\mbm{y}))$

\hspace{2.2cm} $\to\,\,\exists D\forall C(H(C)\to \exists !! \mbm{y}(C\mbm{y}\land D\mbm{y}))$}} 

\vspace{0.2cm}

 If we restrict the set of the classes $k$ considered by Russell in \cite{Russ} to a domain of non-empty $m$-ary predicates with properties expressed by a formula $H(C)$, then Russell's axiom has the form $choice_*^{m}(H)$. In this way, we get the set $choice_*^{(2)}$ of these principles.

\vspace{0.3cm}

\noindent{\bf The set $choice_*^{(2)}$ of all Asser axioms of in HPL}

\nopagebreak

\noindent\fbox{\parbox{11.8cm}{

\vspace{0.1cm} 
$\begin{array}{lll}
choice_*^{(2)}&=_{\rm df}\bigcup_{m\geq 1}choice_*^{m}\\
choice_*^{m}&=_{\rm df}\{choice_*^{m}(H)\mid H \mbox{ is a formula in ${\cal L}_C^{(2)}$} \}\\
\end{array}
$}} 

\vspace{0.2cm} 
For $m\geq 2$, $choice_*^{m}(H)$ is stronger than $AC_*^{1,1}$. 

\section{Outlook: The independence of the second-order WO}\label{SectionIndependWO}

\begin{proposition}\label{choiceInS1} For any $n, m\geq 1$ and any $H$ in ${\cal L}_{\mbmss{x},D}^{(2)}$, there holds 
\[\S_0\models choice_h^{n,m}(H).\]
\end{proposition}

\begin{proposition}\label{WO_unabh_choice} There is a Henkin-Asser structure that is a model of 
\[^{h}ax^{(2)} \cup choice_h^{n,m} \cup \{\neg WO^1\}\]
 for any $n\geq 1$ and $m\geq 1$.
\end{proposition}

\setcounter{satz}{0}

Proposition \ref{WO_unabh_choice} answers a question of Paolo Mancosu and Stewart Shapiro (2018). 

By Proposition \ref{WO_unabh_choice} there is a Henkin-Asser structure that is a model of $\{\neg WO^{1}\}\cup choice_h^{n,m}$. Moreover, if we use ZFC as metatheory and ZFC is consistent, then every standard structure $\S_1^I$ is a model of $WO^{1}$ by \cite{Zerm04} and thus a model of $\{ WO^{1}\}\cup choice_h^{n,m}$. 

\begin{theorem}\label{HPL_ind_3}$WO^1 $ is {\rm HPL}-independent of $ choice_h^{n,m}$.
\end{theorem}

Note, that we also have the following theorems.
\begin{theorem}\label{HPL_ind_1} $choice_*^1$ is {\rm HPL}-independent of $AC^{1,1}$. 
\end{theorem}
\begin{theorem} \label{HPL_ind_2} $choice_h^{1,1}$ is {\rm HPL}-independent of $AC^{1,1}$. \end{theorem}


\begin{thebibliography}{76}
\addcontentsline{toc}{chapter}{\bf Bibliography}

\bibitem{Ack111} {\sc Ackermann,\,\,W.},\,\,\,{\sl Zum Eliminationsproblem der mathematischen Logik}. {\sf Mathema\-tische Annalen} 111, 1935, pp. 61--63.

\bibitem{1a} {\sc Asser,\,\,G.},\,\,\,{\sf Einf\"uhrung in die mathematische Logik, Teil I: Aussagenkalk\"ul}.
 Teubner,\,\,1972.

 \bibitem{1b} ---, ---, {\sf Teil II: Pr\"adikatenkalk\"ul der ersten Stufe}.
 Teubner,\,\,1972.

\bibitem{1c} ---, ---, {\sf Teil III: Pr\"adikatenlogik h\"oherer Stufe}.

\bibitem{Church32} {\sc Church,\,\,A.},\,\,\,{\sl A Set of Postulates for the Foundation of Logic}.\,\,{\sf Annals of Mathematics} 33 (2), 1932, pp. 346--366.

\bibitem{FA22a} {\sc Fraenkel,\,\,A.},\,\,\,{\sl Der Begriff ``definit'' und die Unabh\"angigkeit des Auswahl\-axioms},\,\,\,{\sf Sitzungsberichte der K\"oniglich Preussischen Akademie der Wissenschaften, phys.-math. Kl.}, 1922, pp. 253--257.

\bibitem{FA37} {\sc Fraenkel,\,\,A.},\,\,\,{\sl Ueber eine abgeschwaechte Fassung des Auswahlaxioms},\,\,\,{\sf The Journal of Symbolic Logic} 2 (1), 1937, pp. 1--25.

\bibitem{Frege} {\sc Frege,\,\,G.},\,\,\,{\sf Begriffsschrift} (1879). In M.Wille: Gottlob Frege Begriffsschrift, eine der arithmetischen nachgebildete Formelsprache des reinen Denkens. Springer,\,\,2018.

\bibitem{Henk} {\sc Henkin,\,\,L.},\,\,\,{\sl Completeness in Theory of Types}. {\sf Journal of Symbolic Logic} 14 (3), 1950, pp. 81--91.

 \bibitem{HiAc38} {\sc Hilbert,\,\,D.},\,\,and {\sc W.\,Ackermann},\,\,\,{\sf Grundz\"uge der theoretischen Logik}.  2nd Edition,\,\,	 Springer,\,\,1938. 

 \bibitem{HiAc59} ---, ---. 4th Edition,\,\,	 Springer,\,\,1959. 

 \bibitem{HiBe34} {\sc Hilbert,\,\,D.},\,\,and P.\,Bernays,\,\,\,{\sf Grundlagen der Mathematik I}. In: Die Grundlehren der mathematischen Wissenschaften in Einzeldarstellungen mit besonderer Berücksichtigung der Anwendungsgebiete. Volume XL, Springer, 1934.

\bibitem{Gass84} {\sc Gassner,\,\,C.},\,\,\,{\sf Das Auswahlaxiom im Pr\"adikatenkalk\"ul
zweiter Stufe}. Dissertation, Greifs\-wald, 1984.

 \bibitem{Gass94} {\sc Gassner,\,\,C.},\,\,\,{\sl The Axiom of Choice in Second-Order Predicate Logic}. {\sf Mathematical Logic Quarterly}, 40 (4), 1994, pp. 533--546.

\bibitem{Jech73} {\sc Jech,\,\,T. J.},\,\,\,{\sf The Axiom of Choice} (1973). Dover Publications, 2008.

\bibitem{M38} {\sc Mostowski,\,\,A.},\,\,\,{\sl Über den Begriff einer endlichen Menge}.\,\,{\sf Comptes rendus des s\'eances de la Soci\'et\'e des Sciences et des Lettres de Varsovie},\,\,Classe III,\,\,31, 1938, pp. 13--20.

\bibitem{M39} {\sc Mostowski,\,\,A.},\,\,\,{\sl Über die Unabh\"angigkeit des Wohlordnungssatzes vom Ordnungsprinzip}. {\sf Fundamenta Mathematicae} 32, 1939, pp. 201--252.

\bibitem{Russ} {\sc Russell,\,\,B.} (1906),\,\,\,{\sl On Some Difficulties in the Theory of Transfinite Numbers and Order Types}. {\sf Proc. London Math. Soc.} Volume s2-4,\,\,Issue 1,\,\,1 January 1907 pp. 29--53. 

\bibitem{DudReMa}{\sc Scheid,\,\,H.},\,\,\,{\sf Duden. Rechnen und Mathematik. Das Lexikon für Schule und Praxis}. Dudenverlag,\,\,1994.

\bibitem{Siskind} {\sc Siskind,\,\,B.},\,\,\,{\sc Mancosu,\,\, P.}, and {\sc Shapiro,\,\,S.},\,\,\,{\sl A Note on Choice Principles in Second-Order Logic}. Review of Symbolic Logic 16 (2), 2023, pp. 339--350. 

 \bibitem{Vaeae} {\sc V\"a\"an\"anen,\,\,J.\,A.},\,\,\,{\sl Second-order and higher-order logic}. {\sf Stanford Encyclopedia of Philosophy}. Zalta, E. N. (ed.). Metaphysics Research Lab, Stanford, (Stanford encyclopedia of philosophy), August, 2019.

 \bibitem{Zerm04} {\sc Zermelo,\,\,E.}, {\sl Beweis,\,\,da\ss{} jede Menge wohlgeordnet werden kann} (Aus einem an Herrn Hilbert gerichteten Briefe).
{\sf Mathematische Annalen} 59 (Zeitschriftenheft 4), 1904, pp. 514--516.
 Teubner,\,\,1981.

\end{thebibliography}
\end{document}